[1994/12/01]

\documentclass[12pt]{article}
\usepackage{amscd}
\usepackage{amsthm} 
\usepackage{amsfonts}
\usepackage{amssymb}
\usepackage{amsmath}

\newcommand{\field}[1]{\ensuremath{\mathbb{F}_{#1}}} 

\begin{document}

\title{Yet More Projective Curves Over $F_2$}
\author{ Chris Lomont, Purdue University
\date{Dec 30, 2000}
\thanks{clomont@math.purdue.edu}}
\maketitle

\begin{abstract} All binary plane curves of degree less than 7 are
examined for curves with a large number of \field{q} rational
points on their smooth model, for $q=2^m, m = 3,4,...,11$. Previous results are improved,
and many new curves are found meeting or close to Serre's, Lauter's, and
Ihara's upper bounds.
\end{abstract}

\section{Introduction}

All absolutely irreducible homogeneous polynomials $f \in
\field{2}[x,y,z]$ of degree less than 7 are examined for those with a
large number of \field{q} rational points, extending the results in
\cite{Moreno}.

Let \field{q} denote the finite field of $q=2^m$ elements. Let $f
\in \field{q}[x,y,z]$ be an absolutely irreducible homogeneous
polynomial; $f$ defines a projective plane curve $C$. Let
$\tilde{C}$ be the smooth model, and $g$ its genus. Some bounds on the genus can be deduced from knowing information about the number $N_q$ of \field{q} rational points. Serre's bound
\cite{Serre} on the number $N_q$ of \field{q}-rational points of $\tilde{C}$ is $$|N_q - (q+1)|\leq g \lfloor 2\sqrt{q}
\rfloor$$ where $\lfloor \alpha \rfloor$ is the integral part of
$\alpha$. This gives $$\frac{|N_q - (q+1)|}{\lfloor 2\sqrt{q}\rfloor} \leq g$$ Also if there is an integer $g_0$ such that $$q+1+g_0\lfloor 2 \sqrt{q} \rfloor < N_q$$ then $g_0<g$. If the number of singularities is $r$, and the degree of the plane curve is d, then $$g\leq \frac{(d-1)(d-2)}{2}-r$$ This
gives some upper and lower bound on the possible genus, which speeds
up the computation by removing curves with low numbers of rational points early in the computation. To get an
estimate of the total number of points possible on the smooth
model $\tilde{C}$ resulting from blowing up singularities the
following simple estimate was used.

\newtheorem{thrm1}{Theorem}
\newtheorem{def1}{Definition}

\begin{thrm1} Let $C\subseteq \mathbb{P}^2$
be a plane curve of degree d with singularities $P_1, P_2,
...,P_r$, with multiplicities $m_1,m_2,...,m_r$, for $r\geq 2$. Then
$\sum^r_{i=1}m_i \leq \lfloor \frac{d}{2} \rfloor r + 1$ if d is
odd, and $\sum^r_{i=1}m_i \leq \lfloor \frac{d}{2} \rfloor r$ if d
is even.
\end{thrm1}

So the number of points obtained from blowing up singularities is bounded above by $\lfloor \frac{d}{2} \rfloor r + 1$.

\begin{proof}
By Bezout's theorem, a line through any 2 singularities $P_i$ and
$P_j$ implies $m_i + m_j \leq d$. Thus at most one singularity can
have multiplicity $> d/2$, and the result follows.
\end{proof}

More details on resolution of curve singularities can be found in
\cite{Hartshorne} and \cite{Walker}. The number of \field{q} smooth points on a plane curve can be computed directly, as can lower bounds on the number of singularities, and then the above inequalities can be used to obtain bounds on the genus and on $N_q$.

To test for absolute
irreducibility the following was used \cite{Ragot}:

\begin{def1}
Let $k$ be a field. The polynomial $f\in k[x_1,x_2,...,x_n]$ has a
{\bf simple solution} at a point $P\in k^n$ if $f\in I(P)\diagdown
I(P)^2$, with $I(P)$ being the ideal of polynomials vanishing at
$P$.\end{def1}

\begin{thrm1} If $f\in k[x_1,x_2,...,x_n]$ is irreducible over the
perfect field $k$ and has a simple solution in $k^n$, then $f$ is
absolutely irreducible.\end{thrm1}
\begin{proof} Since $f$ is irreducible over $k$, its absolutely
irreducible factors are conjugate over $k$. If $P\in k^n$ is a
root of one of the factors, it must be a root of the others. But
$P$ only vanishes to order 1, thus $f$ has one factor, and is
absolutely irreducible. \end{proof}

The layout of this paper follows \cite{Moreno}.

\section{Computation} (or how to store 268 million polynomials in
32 megabytes of RAM)

All homogeneous polynomials of degree $\leq 5$ in $\field{2}[x,y,z]$ were examined
over the fields \field{q}, for $q=3,4,...,11$. Due to the the
length of the computation, degree 6 homogeneous polynomials were examined only
for $q=3,4,...,9$.

The most time consuming part was counting the number of
\field{q}-rational points on each plane curve. This was done with
a C program using exhaustive search. Several ideas were used to reduce the
complexity at each stage. Degree 6 computations will be described;
the other degrees are similar. The code was checked for correctness by comparing the
degree $\leq$ 5 results with \cite{Moreno}, and in the process a few curves were found that were previously overlooked.

First, each homogeneous polynomial can be represented uniquely by
a 32 bit integer, using each bit to signify the presence of a
certain monomial in the polynomial. In degree 6, there are
$\binom{6+2}{2}=28$ different monomials of the form $x^i y^j z^k$
with $i+j+k=6$ and $0\leq i,j,k$. So each bit from 0 to 27 denotes
the presence of a monomial, and the mapping $\alpha : \{f\in
\field{2}[x,y,z] | f homogeneous, deg(f)=6\} \rightarrow
\{1,2,...,2^{28}-1 \}$ thus defined is a bijection. So each
homogeneous polynomial $f\in \field{2}[x,y,z]$ of degree 6
corresponds to a unique integer between $1$ and $2^{28}-1$.

To reduce the number of polynomials searched, equivalent ones
under the action of $GL_3(\field{2})$ were removed. To fit the
entire degree 6 computation in memory, a bit table of 32 megabytes
of RAM was used, with the position of each bit representing the
number of a polynomial. All bits were set to 1, denoting the
polynomial is still in the search space, and then the orbit of
each polynomial under $GL_3(\field{2})$ was removed from the bit
table, trimming the 268 million polynomials down to 1.6 million.
This showed the four polynomials in \cite{Moreno} of degree 4,
genus 3, with 113 rational points over \field{64}, are 2 distinct
polynomials mod $GL_3(\field{2})$.

For speed reasons finding solutions was done by table lookup, so
in each orbit the polynomial needing the fewest number of lookups
was selected. By choosing the representative with the fewest
number of lookups as opposed to the polynomial with the lowest
value of $\alpha(f)$ defined above, 12 million lookups were
removed from polynomials of degree 6, resulting in over 3 trillion
operations removed during the rational point counting. Also,
clearly reducible polynomials, such as those with all even
exponents or divisible by a variable, are removed at this point.
At each stage data is saved to prevent having to rerun any step.

After the C program computes all the \field{q}-rational points,
the points are tested for singularities. The computation up to
this point took about 80 hours of computer time on a Pentium III
800 MHZ. Using the bounds above on the genus and possible ranges
for number of \field{q} rational points on the smooth model, the
program searches all curves for those with a large number of
possible \field{q} rational points for each genus and field
combination, and all such curves are written out to be examined.
If the genus of one of these curves is not unique using the
bounds, the program KANT \cite{Kant} is used to compute the genus,
and this data is incorporated into the C program, and another pass
is run. Due to the large number of degree 6 curves, and the length
of time to compute the genus of all of them, not all degree 6
curves of genus $\leq 5$ were identified. All curves of degree 6,
genus $\geq 6$ were identified. The C program also found simple
points over \field{2} to apply the theorem above, and then Maple V
\cite{Maple} was used to test for irreducibility since it has
multivariable factoring algorithms over finite fields. For 12
curves of degree 6, there were no simple \field{2} points, so
\field{4} simple points were used. For 2 of these curves there
were no such simple \field{4} points, so \field{8} simple points
were used. This turned out to suffice to check absolute
irreducibility of all polynomials in this paper. The C program
also found the singularity types of the \field{2} singularities
for visual inspection to see if there are clearly more rational
points on the smooth model. The package of \cite{Hache} was not
available to do a more detailed singularity analysis, thus some of
the bounds below may be improved by looking for rational points
over a wider class of singularities than the \field{2}
singularities considered here.

\section{Computational Results}
For each field and genus combination all polynomials are listed
that result in the maximal known number of rational points on the
smooth model of the curve. For fields \field{q}, $q=3,4,..,11$ all
binary homogeneous polynomials of degree $\leq 5$ were searched.
For $q=3,4,...9$ the search was extended to include all degree 6
binary homogeneous polynomials. Let $N_q(g)$ denote the maximal
number of smooth points on all curves of genus $g$ over \field{q}.
For genus and field combinations not listed here, see
\cite{Moreno}.

\subsection{\field{8}}
A new curve of genus 3 and the maximal number of smooth points, 24, is
$$ x^6 + x^5\,y + x\,y^5 + y^6 +
  \left( x^5 + x^2\,y^3 + y^5 \right) \,z + x^3\,y\,z^2 +
  x\,y^2\,z^3 + \left( x + y \right) \,z^5 + z^6$$
A genus 5 curve with 28 planar smooth points is $$x^6 + x^5\,y + x^3\,y^3 + y^6 + y^5\,z +
y^4\,z^2 + \left( x^3 + x\,y^2 + y^3 \right) \,z^3 +
  \left( x^2 + x\,y \right) \,z^4 + x\,z^5$$
A genus 6 curve, with 33 planar smooth points is $$x^4\,y^2 +
x^3\,y^3 + x\,y^5 + x\,y^4\,z + y^4\,z^2 +  \left( x^2 + y^2
\right) \,z^4 + y\,z^5$$ A curve of genus 7 with 33 smooth planar
points is
  $$x^6 + x^5\,y + x^4\,y^2 + x^3\,y^3 + y^6 + y^5\,z +
  \left( x^2\,y^2 + x\,y^3 + y^4 \right) \,z^2 +
  \left( x^2\,y + x\,y^2 \right) \,z^3 +
  \left( x^2 + x\,y + y^2 \right) \,z^4$$
A genus 8 curve with 33 smooth planar points is $$x^5\,y
+ x^2\,y^4 + x\,y^5 + y^6 +
  \left( x^4\,y + x^2\,y^3 \right) \,z + x^3\,y\,z^2 +
  \left( x^3 + x\,y^2 \right) \,z^3 + x\,y\,z^4 + y\,z^5$$
Two curves of genus 9, each with 33 smooth planar points:
$$f_1=x^5\,y + x^4\,y^2 + x^2\,y^4 +
  \left( x^3\,y^2 + x^2\,y^3 \right) \,z + x^4\,z^2 +
  \left( x\,y^2 + y^3 \right) \,z^3 + x^2\,z^4 + x\,z^5 +
  z^6$$ $$f_2=x^6 + x^4\,y^2 + x^3\,y^3 + x^2\,y^4 + x^3\,y^2\,z +
  \left( x^4 + x^3\,y + x\,y^3 \right) \,z^2 + y^3\,z^3 +
  \left( x + y \right) \,z^5 + z^6$$
Five curves of genus 10 with 35 smooth points in the plane:
  $$f_1 = x^5\,y + x^2\,y^4 + y^6 +
  \left( x^3\,y^2 + x\,y^4 + y^5 \right) \,z + y^4\,z^2 +
  \left( x^2\,y + x\,y^2 \right) \,z^3 +
  \left( x^2 + y^2 \right) \,z^4 + x\,z^5$$ $$f_2 = x^5\,y + x^4\,y^2 +
  x^3\,y^3 + x\,y^5 + y^6 + x^2\,y^3\,z + \left( x^4 + x^2\,y^2 + x\,y^3 \right) \,
   z^2 + x^3\,z^3 + y^2\,z^4 + x\,z^5 + z^6$$ $$f_3 = x^5\,y + x^4\,y^2 + x^2\,y^4 +
  \left( x^4\,y + y^5 \right) \,z +
  \left( x^4 + x\,y^3 \right) \,z^2 + x^3\,z^3 +
  \left( x^2 + y^2 \right) \,z^4 + y\,z^5 + z^6 $$ $$f_4 = x^5\,y + x^3\,y^3 + x^2\,y^4 +
  \left( x^5 + x^2\,y^3 + y^5 \right) \,z +
  \left( x^2\,y + y^3 \right) \,z^3 + x^2\,z^4 +
  \left( x + y \right) \,z^5 + z^6 $$ $$f_5 = x^4\,y^2 + x^2\,y^4 + x\,y^5 +
  x^5\,z + x\,y^3\,z^2 + \left( x^3 + x^2\,y \right) \,z^3 +
  \left( x^2 + y^2 \right) \,z^4 + x\,z^5 $$

\subsection{\field{16}}
One genus 6 curve, a Hermitian curve, with the maximal number of
smooth points, 65, was found: $$x^5 + y^5 + z^5$$ Two genus 7
curves each with 57 smooth points in the plane: \begin{eqnarray*}
f_1&=&x^4\,y^2 + x\,y^5 + y^6 +
  \left( x^2\,y^3 + x\,y^4 + y^5 \right) \,z +
  \left( x^2 + x\,y \right) \,z^4 + y\,z^5 \\  f_2&=&x^4\,y^2 + x^2\,y^4 + y^6 +
  \left( x^2\,y^3 + x\,y^4 \right) \,z +
  \left( x^2 + x\,y \right) \,z^4 + y\,z^5\end{eqnarray*}
A curve with genus 8 with 57 smooth plane points is $$x^6 +
x^3\,y^3 + x^2\,y^4 + x^4\,y\,z + x^2\,y^2\,z^2 + \left( x^3 +
x^2\,y \right) \,z^3 +
  \left( x + y \right) \,z^5$$
There are two curves of genus 9 with 57 smooth plane points, each
receiving two points from blowups: $f_1$ from the singularity
$(1:1:1)$ of type $u^2+u v + v^2$ which splits over \field{16},
and $f_2$ from the singularity $(0:1:1)$ of type $u v$. Thus
$N_{16}(9)\geq 59$.
\begin{eqnarray*} f_1&=&x^5\,y + x^3\,y^3 + x\,y^5 +
  \left( x^5 + y^5 \right) \,z + x^2\,y^2\,z^2 +
  \left( x^3 + y^3 \right) \,z^3 +
  \left( x + y \right) \,z^5 \\ f_2 &=& x^6 + x^5\,y + x^2\,y^4 + y^5\,z + x^2\,y^2\,z^2 +
  x\,y^2\,z^3 + \left( x^2 + x\,y \right) \,z^4 + y\,z^5\end{eqnarray*}
The two curves of genus 10 each with 59 plane smooth points are:
  $$f_1 = x^5\,y + y^6 + \left( x^2\,y^3 + y^5 \right) \,z +
  \left( x^4 + x^3\,y + x\,y^3 \right) \,z^2 +
  x\,y^2\,z^3 + \left( x + y \right) \,z^5 + z^6$$
  $$f_2 = x^5\,y + y^6 + \left( x^4\,y + x\,y^4 + y^5 \right) \,z +
  \left( x^4 + x\,y^3 \right) \,z^2 +
  \left( x^3 + x\,y^2 \right) \,z^3 + y^2\,z^4 + z^6$$

\subsection{\field{32}}
Three curves with genus 4 and 71 smooth points on the plane curve
are: \begin{eqnarray*} f_1 &=& x^4\,y + x\,y^4 + y^5 + x\,y^3\,z +
  \left( x\,y^2 + y^3 \right) \,z^2 + x^2\,z^3 + x\,z^4 +
  z^5 \\ f_2&=& x^6 + x^3\,y^3 + y^6 + \left( x^4\,y + y^5 \right) \,z +
  \left( x^3\,y + x^2\,y^2 \right) \,z^2 +
  \left( x^3 + x^2\,y + y^3 \right) \,z^3 + \\&& x^2\,z^4 +
  y\,z^5 + z^6 \\ f_3 &=& x^6 + x^3\,y^3 + y^6 + \left( x^5 + x^3\,y^2 +
     x^2\,y^3 \right) \,z + y^4\,z^2 +
  \left( x^3 + y^3 \right) \,z^3 + y^2\,z^4 + x\,z^5 + z^6
  \end{eqnarray*}
A curve with 82 smooth points in the plane and genus 5 is
  $$x^6 + x^3\,y^3 + x^2\,y^4 + y^6 + x^5\,z + x^3\,y\,z^2 +
  \left( x^3 + x\,y^2 + y^3 \right) \,z^3 + x^2\,z^4 +
  y\,z^5$$
A genus 6 curve with 82 planar smooth points and 2 points above
the singularity $(1:0:1)$ of type $u v$ (thus $N_{32}(6)\geq 84$)
is $$x^6 + y^6 + \left( x^4\,y + x^3\,y^2 + x\,y^4 \right) \,
   z + x\,y^2\,z^3 + \left( x^2 + x\,y + y^2 \right) \,z^4$$
Two genus 7 curves each with 92 planar smooth points are
\begin{eqnarray*}f_1&=&x^3\,y^3 + y^6 + \left( x^5 + x^3\,y^2 \right) \,z +
  \left( x^4 + y^4 \right) \,z^2 +
  \left( x^3 + y^3 \right) \,z^3 + y^2\,z^4 + x\,z^5 + z^6\\
   f_2&=&x^6 + y^6 + \left( x^5 + y^5 \right) \,z + y^4\,z^2 +
  \left( x^3 + y^3 \right) \,z^3 + x^2\,z^4 +
  \left( x + y \right) \,z^5\end{eqnarray*}
A curve with 93 planar smooth points, genus 8, is $$x\,y^5 + y^6 + \left( x^5 + x^4\,y
\right) \,z + y^4\,z^2 + \left( x^3 + y^3 \right) \,z^3 + y^2\,z^4
+ y\,z^5$$
A genus 9 curve with 93 smooth planar points:  $$x^4\,y^2 +
x^3\,y^3 + \left( x^5 + x^3\,y^2 + x\,y^4 +
     y^5 \right) \,z + x^2\,y^2\,z^2 +
  \left( x^3 + y^3 \right) \,z^3 + x^2\,z^4 + z^6$$
Genus 10 with 103 smooth planar points:
  $$x^6 + x^3\,y^3 + x\,y^5 +
  \left( x^2\,y^2 + x\,y^3 \right) \,z^2 +
  \left( x^3 + x\,y^2 + y^3 \right) \,z^3 + x\,y\,z^4 +
  \left( x + y \right) \,z^5$$

\subsection{\field{64}}
One curve had genus 4 and 118 smooth planar points: $$x^3\,y^2 +
y^5 + y^4\,z + y^2\,z^3 + z^5$$
Two curves of genus 6 had 160 smooth planar
points (which is one less than the bound of 161): \begin{eqnarray*} f_1 &=& x^4\,y^2 + x^2\,y^4 + x\,y^5 + y^5\,z +
y^3\,z^3 +
  y\,z^5 + z^6 \\ f_2&=&x^6 + x^5\,z + \left( x^4 + y^4 \right) \,z^2 + x^3\,z^3 +
  y^2\,z^4 + y\,z^5\end{eqnarray*}
Genus 7, 153 planar smooth points: $$x^2\,y^4 + x\,y^5 + y^6 +
  \left( x^3\,y^2 + x\,y^4 + y^5 \right) \,z +
  x\,y^3\,z^2 + x^2\,z^4 + x\,z^5 + z^6$$
Three curves had genus 8 and 159 plane smooth points, the last two
of which have no rational points over \field{2}:
\begin{eqnarray*} f_1 &=& x^3\,y^3 + y^6 + \left( x^4\,y + x\,y^4
\right) \,z +
  \left( x^3 + y^3 \right) \,z^3 + x\,y\,z^4 \\ f_2 &=& x^6 + x^5\,y + x^3\,y^3 + x\,y^5 + y^6 +
  \left( x^3\,y^2 + y^5 \right) \,z + x^3\,y\,z^2 +
  \left( x^2\,y + y^3 \right) \,z^3 + \\&& y^2\,z^4 + x\,z^5 +
  z^6 \\ f_3 &=& x^6 + x^4\,y^2 + x^3\,y^3 + x^2\,y^4 + y^6 +
  \left( x^4 + x^2\,y^2 + y^4 \right) \,z^2 +
  \left( x^3 + y^3 \right) \,z^3 + \\&&
  \left( x^2 + y^2 \right) \,z^4 + z^6\end{eqnarray*}
There are 166 plane smooth points on this curve of genus 9: $$x^6
+ x^3\,y^3 + \left( x^4\,y + x^2\,y^3 \right) \,z +
  \left( x^3\,y + x\,y^3 + y^4 \right) \,z^2 + x^2\,z^4 +
  y\,z^5$$
Four genus 10 curves each had 171 points on their smooth model:
\begin{eqnarray*} f_1 &=& x^6 + y^6 + \left( x^4\,y + x^2\,y^3 + x\,y^4 \right)
\,
   z + x^3\,y\,z^2 + x\,y^2\,z^3 + x\,y\,z^4 + z^6 \\ f_2 &=& x^6 + x^5\,y + x^4\,y^2 + x^3\,y^3 + x^2\,y^4 + x\,y^5 +
  y^6 + \left( x^4\,y + x\,y^4 \right) \,z +
  \left( x^2\,y + x\,y^2 \right) \,z^3 + z^6\\ f_3 &=& x^6 + x^3\,y^3 + y^6 + \left( x^4\,y + x\,y^4 \right) \,
   z + \left( x^3 + y^3 \right) \,z^3 + z^6 \\f_4 &=& x^6 + x^3\,y^3 + x\,y^5 + x^3\,y^2\,z +
  \left( x^4 + x^3\,y + y^4 \right) \,z^2 + y^2\,z^4 +
  x\,z^5 + z^6 \end{eqnarray*}

\subsection{\field{128}}
There is one degree 6 plane curve with a genus 3 smooth model, with 183 smooth plane points, and
another point coming from the singularity $(0:0:1)$ of type $(u+v)(u^2+u v +
v^2)$, which matches \cite{Moreno}. The curve is $$x^6 + x^5\,y + x^4\,y^2 + x^3\,y^3 +
x^2\,y^4 +
  \left( x^5 + x^4\,y \right) \,z + y^4\,z^2 +
  \left( x^3 + y^3 \right) \,z^3$$
A curve of genus 4 with 215 planar smooth points (2 less than the maximum possible) is $$x^2\,y^3 +
x\,y^4 + x^4\,z + x\,y^2\,z^2 + x\,y\,z^3 +
  \left( x + y \right) \,z^4$$
There are two curves of genus 6 with 240 planar smooth points,
receiving 3 points each from singularities. $f_1$ has type $u
v(u+v)$ at $(0:1:0)$ and $f_2$ has type $(u+v)(u^2+u v+v^2)$ at
$(0:1:0)$ and type $u v$ at $(1:0:0)$. Thus $N_{128}(6)\geq 243$.
\begin{eqnarray*} f_1 &=& x^4\,y^2 + \left( x^5 + x^4\,y + x^2\,y^3 \right) \,z +
  \left( x^2\,y^2 + x\,y^3 \right) \,z^2 +
  \left( x^2 + x\,y + y^2 \right) \,z^4 \\ f_2 &=& x^3\,y^3 + x^4\,y\,z + \left( x^4 + x^3\,y \right) \,z^2 +
  \left( x^3 + x^2\,y + y^3 \right) \,z^3 + z^6\end{eqnarray*}
Two genus 7 curves with 248 smooth planar points:
\begin{eqnarray*} f_1 &=& x^3\,y^3 + x\,y^5 + y^6 + x^3\,y^2\,z + y^4\,z^2 +
  x^3\,z^3 + \left( x^2 + y^2 \right) \,z^4 +
  \left( x + y \right) \,z^5 \\ f_2 &=& x^5\,y + x^4\,y^2 + x^2\,y^4 +
  \left( x^3\,y^2 + x^2\,y^3 \right) \,z + x^4\,z^2 +
  x\,y^2\,z^3 + z^6 \end{eqnarray*}
A curve with 266 planar smooth points, genus 8, and no \field{2}
rational points is $$x^6 + x^3\,y^3 + x^2\,y^4 + x\,y^5 + y^6 +
  \left( x^5 + y^5 \right) \,z +
  \left( x^2\,y^2 + y^4 \right) \,z^2 +
  \left( x^3 + y^3 \right) \,z^3 + x\,z^5 + z^6$$
There are 269 smooth plane points on the curves of genus 9 given
by \begin{eqnarray*} f_1 &=& x^4\,y^2 + x\,y^5 + \left( x^4 + y^4
\right) \,z^2 +
  \left( x^3 + y^3 \right) \,z^3 + x\,y\,z^4 + x\,z^5 + z^6 \\ f_2
  &=& x^6 + x^3\,y^3 + x^2\,y^4 + y^6 +
  \left( x\,y^4 + y^5 \right) \,z + x^2\,y^2\,z^2 +
  \left( x^2 + x\,y + y^2 \right) \,z^4 \end{eqnarray*}
The smooth curve of genus 10 with 276 \field{128} rational points
is $$x^6 + y^6 + x^2\,y^3\,z +
  \left( x^4 + x^3\,y + y^4 \right) \,z^2 + x^3\,z^3 +
  x^2\,z^4 + x\,z^5$$

\subsection{\field{256}}
A genus 3 curve not listed in \cite{Moreno} with 350 smooth planar
points is given by $$x^5 + x\,y^4 + y^5 + \left( x^2\,y^2 + y^4
\right) \,z +
  \left( x^2\,y + x\,y^2 \right) \,z^2 + x\,z^4 + z^5$$
A curve with 399
smooth plane points and genus 5 is $$x^6 + x^4\,y^2 + x^5\,z +
  \left( x^2\,y^2 + y^4 \right) \,z^2 +
  \left( x^2\,y + x\,y^2 \right) \,z^3 + x^2\,z^4 + y\,z^5$$
A genus 6 curve with 416 smooth plane points is $$x^4\,y +
x^3\,y^2 + y^4\,z +
  \left( x^3 + y^3 \right) \,z^2 +
  \left( x^2 + x\,y \right) \,z^3 + z^5$$
One point from the singularity $(1:0:0)$ of type $u v^2$ is
added to the 442 smooth plane points on a curve of genus 7 given
by $$x^3\,y^3 + x^2\,y^4 + y^5\,z + x^3\,y\,z^2 +
  \left( x\,y^2 + y^3 \right) \,z^3 + y^2\,z^4 + z^6$$
A curve of genus 8 with one point less than the Serre bound has
512 smooth plane points and is given by $$x^4\,y^2 + y^5\,z +
x\,z^5$$ Two curves of genus 9, each with 474 smooth points and 2
points from singularities of type $u^2+u v + v^2$, which factor
over \field{256}, at points $(0:1:1)$ and $(1:1:0)$ respectively
(so $N_{256}(9) \geq 476$) are \begin{eqnarray*}f_1&=&x^5\,y +
x^3\,y^3 + x^2\,y^4 + x\,y^5 + y^4\,z^2 + x^3\,z^3 + y^2\,z^4 +
x\,z^5 \\f_2&=&x^6 + y^6 + \left( x^5 + y^5 \right) \,z + x^4\,z^2
+ \left( x^2\,y + x\,y^2 \right) \,z^3 + x\,z^5 +
z^6\end{eqnarray*} Two smooth curves of genus 10 have 537 smooth
plane points: \begin{eqnarray*} f_1&=&x^6 + x\,y^5 + x^4\,y\,z +
x^2\,y^2\,z^2 + y^3\,z^3 +
  x\,z^5 \\f_2&=&x^6 + x^5\,y + x^3\,y^3 + x\,y^5 + y^6 +
  \left( x^5 + y^5 \right) \,z + x^2\,y^2\,z^2 +
  x\,y\,z^4 + z^6\end{eqnarray*}
\subsection{\field{512}}
Four curves overlooked in \cite{Moreno} of genus 4 have 663 plane smooth points. They are
\begin{eqnarray*}f_1 &=& x^4\,y + x\,y^4 + \left( x^3\,y + y^4 \right) \,z +
  \left( x\,y + y^2 \right) \,z^3 + z^5 \\ f_2&=&x^4\,y + x\,y^4 + y^5 +
  \left( x\,y^3 + y^4 \right) \,z + \left( x\,y^2 + y^3 \right) \,z^2 +
  x^2\,z^3 + x\,z^4 \\ f_3&=&x^4\,y + x^2\,y^3 + y^5 +
  \left( x^2\,y^2 + x\,y^3 + y^4 \right) \,z +  \left( x^3 + y^3 \right)
  \,z^2 + z^5 \\ f_4&=&x^5 + y^5 + \left( x^4 + x^3\,y + y^4 \right) \,z +
  \left( x\,y^2 + y^3 \right) \,z^2 + z^5\end{eqnarray*}
A genus 6 curve with 766 smooth plane points and one more point
from the singularity $(1:0:0)$ of type $(u+v)(u^2+u v+v^2)$ (so
$N_{512}(6)\geq 767$) is $$x^3\,y^3 + y^6 + \left( x\,y^4 + y^5
\right) \,z +
  \left( x^2\,y^2 + y^4 \right) \,z^2 + x^3\,z^3 +
  x\,y\,z^4 + \left( x + y \right) \,z^5$$
There are 786 smooth plane points and 1 point from the
singularity $(1:0:0)$ of type $u v^2$ on the genus 7 curve
$$x^2\,y^4 + y^6 + x^3\,y^2\,z +
  \left( x^3 + x\,y^2 \right) \,z^3 + x\,y\,z^4 + y\,z^5$$
A curve of genus 8 with 813 plane smooth points is $$x^2\,y^4 +
y^6 + \left( x^5 + x^2\,y^3 \right) \,z +
  \left( x^3\,y + x\,y^3 + y^4 \right) \,z^2 +
  \left( x^3 + x^2\,y \right) \,z^3 + x\,z^5$$
A genus 9 curve with 837 smooth plane points is $$x^6 + x^4\,y^2 +
\left( x\,y^4 + y^5 \right) \,z +
  x^2\,y^2\,z^2 + \left( x^2\,y + x\,y^2 \right) \,z^3 +
  x\,y\,z^4 + x\,z^5 + z^6$$
A smooth genus 10 plane curve with 845 plane points is $$x^5\,y +
x^4\,y^2 + x^2\,y^4 + y^6 +
  \left( x^2\,y^3 + y^5 \right) \,z +
  \left( x^3\,y + y^4 \right) \,z^2 + x^2\,z^4 +
  \left( x + y \right) \,z^5$$

\subsection{\field{1024}}
A genus 3 curve with 1211 smooth plane points is $$x^3\,y + y^3\,z
+ y\,z^3 + z^4$$ Three genus 4 curves have 1273 smooth plane
points:
\begin{eqnarray*} f_1 &=& x^3\,y^2 + y^5 + x^2\,y\,z^2 + y^2\,z^3 + x\,z^4\\
f_2 &=& x^4\,y + x^2\,y^3 + y^5 +
  \left( x^2\,y^2 + y^4 \right) \,z +
  \left( x^3 + y^3 \right) \,z^2 + x\,z^4 \\
  f_3 &=& x^5 + x^3\,y^2 + x^2\,y^3 + y^5 + y^4\,z + x\,y^2\,z^2 +
  x^2\,z^3\end{eqnarray*}
A curve with 1343 smooth plane points, genus 5, and 2 points
coming from the singularity $(1:0:0)$ of type $u v$ (and thus attaining the maximum possible 1345) is
$$x^3\,y^2 + y^5 + x^3\,y\,z + y^3\,z^2 + z^5$$
A genus 6 curve
with 1383 smooth plane points is $$x^4\,y + x\,y^4 + y^5 +
x^2\,y^2\,z + x\,y\,z^3 + z^5$$

\subsection{\field{2048}}
Two genus 3 curves with 2293 smooth plane points, and one more
coming from the singularity $(0:1:0)$ of type $(u+v)(u^2+u v+v^2)$
on each curve are $$f_1=x^4\,y + x^3\,y^2 + x^3\,y\,z + y^2\,z^3 +
  \left( x + y \right) \,z^4$$ $$f_2=x^4\,y + x^3\,y^2 + x^3\,y\,z + x^3\,z^2 + y^2\,z^3 +
  y\,z^4$$
Three curves with 2380 smooth plane points and genus 4 are
\begin{eqnarray*}f_1&=&x^5 + y^5 + x^3\,y\,z + y^2\,z^3 + x\,z^4 \\ f_2&=&x^5 +
x^4\,y + \left( x\,y^3 + y^4 \right) \,z + x^3\,z^2 + \left( x^2 +
y^2 \right) \,z^3 + z^5 \\ f_3&=&x^5 + x^2\,y^3 + x\,y^4 + x^4\,z
+ x^3\,z^2 +  \left( x^2 + x\,y + y^2 \right) \,z^3 + x\,z^4 + z^5
  \end{eqnarray*}
A genus 5 curve with 2422 smooth plane points is $$x^4\,y +
x^3\,y^2 + y^5 + x\,y^2\,z^2 +
  \left( x^2 + x\,y + y^2 \right) \,z^3$$
Finally, a genus 6 curve with 2556 planar smooth points is
$$x^4\,y + x^2\,y^3 + x\,y^4 + y^5 +
  \left( x^3\,y + y^4 \right) \,z + x^2\,z^3 + y\,z^4$$

\subsection{Tallies}

The columns headed ``bound" give the Serre \cite{Serre} bound,
unless marked \cite{Ihara} as Ihara or \cite{Lauter} as Lauter.
The columns headed ``best" give the lower bounds for $N_q(g)$
found above; bounds marked \cite{Moreno} and \cite{Serre} are from
those previous papers. Note in particular the reduction from the
Serre \cite{Serre} upper bounds using \cite{Ihara} and
\cite{Lauter} has made several known curves closer to or already
optimal.

\begin{tabular}{|c|c|c|c|c|c|c|c|c|} \hline
$F_q$ & best 3 & bound  3 & best 4 & bound 4 & best 5 & bound 5 &
best 6 & bound 6
\\ \hline

8 & 24\cite{Serre} & 24 & 25\cite{Moreno} & 28\cite{Lauter} & 28 & 32
\cite{Ihara} & 33 & 35 \cite{Lauter}
\\ \hline

16 & 38\cite{Serre} & 41 & 45\cite{Moreno} & 46\cite{Lauter} & \cite{Moreno} & 54\cite{Lauter} & 65 & 65 \\ \hline

32 & 63\cite{Moreno} & 65\cite{Lauter} & 71 & 76\cite{Lauter} & 82 & 87\cite{Lauter} & 84 & 98\cite{Lauter} \\ \hline

64 & 113\cite{Moreno} & 113 & 118 & 129 & 130\cite{Moreno} & 145 & 160
& 161
\\ \hline

128 & 184\cite{Moreno} & 195 & 215 & 217 & 227\cite{Moreno} & 239 & 243 & 261
\\ \hline

256 & 350\cite{Moreno} & 353 & 381\cite{Moreno} & 385 & 399 & 417 & 416 & 449
\\ \hline

512 & 640\cite{Moreno} & 648 & 663 & 693 & 724\cite{Moreno} & 738
& 767 & 783 \\ \hline

1024 & 1211 & 1217 & 1273 & 1281 & 1345 & 1345 & 1383 & 1409 \\
\hline

2048 & 2294 & 2319 & 2380 & 2409 & 2422 & 2499 & 2556 & 2589 \\
\hline

\hline

$F_q$ & best 7 & bound  7 & best 8 & bound 8 & best 9 & bound 9 &
best 10 & bound 10
\\ \hline

8 & 33 & 39 \cite{Ihara} & 33 & 43\cite{Ihara} &33 &
47\cite{Ihara} &35 & 50\cite{Ihara} \\ \hline

16 & 57 & 70\cite{Ihara} &57 &76\cite{Ihara} &59 &81\cite{Ihara}
&59 &87\cite{Ihara}
\\ \hline

32 & 92 & 110 & 93 &121 & 93 &132  & 103 & 143\\ \hline

64 & 153 &177 & 159 &193 & 166 &209  & 171  & 225\\ \hline

128 & 248 & 283 & 266 &305 & 269 & 327 & 276 & 349 \\ \hline

256 & 443 & 481 & 512 & 513 & 476 & 545 & 537 & 577 \\ \hline

512 & 787 & 828 & 813 & 873  & 837 & 918  & 845 & 963 \\ \hline

\end{tabular}

\section{Conclusions}
Previously known results are \cite{Moreno} and \cite{Serre}. The
bounds in this paper could possibly be strengthened by by
analyzing the singularities in more detail, resulting in more
known \field{q} rational points on the smooth models of the
curves. Also all genus $\leq 5$ curves from the degree 6
polynomials were not identified, although perhaps all the ones
with a large number of \field{q} rational points were. More work
could be done to compute exact parameters for these curves.

\section{Comments}
The techniques used here make a search over degree 7 plane curves
feasible on a supercomputer, and quite possibly on a home PC. The
desingularized curves can be used to construct algebraic-geometric
Goppa codes \cite{Pretzel}, \cite{Tsfasman}. For example, using the genus 5 curve over \field{1024} with 1345 rational points, linear codes with parameters $[n,k-4,n-k]$ can be constructed for $10\leq n \leq 1344$ and $8<k<n$ over \field{1024}. Similarly, using the genus 6 \field{64} curve with 160 points, $[n,k-5, n-k]$-linear codes can be constructed for $12\leq n \leq 159$ and $10< k <n$, and the \field{256} curve of genus 8 with 512 rational points gives $[n,k-7,n-k]$-linear codes for $16\leq n\leq 511$ and $14<k<n$ (see, for example, \cite{Pretzel}).

\end{document}